\documentclass[13pt]{article}
\usepackage{amsmath}
\usepackage{amscd}
\usepackage{amsfonts}
\usepackage{amssymb}
\usepackage{amsthm}
\usepackage[dvips]{color}
\usepackage{epsfig}
\usepackage{graphics}
\usepackage{slashbox}
\usepackage[all,2cell,dvips]{xy} \UseAllTwocells \SilentMatrices
 \newtheorem{lem}{Lemma}[section]
 \newtheorem{cor}[lem]{Corollary}
 \newtheorem{thm}[lem]{Theorem}
 \newtheorem{prop}[lem]{Proposition}
 \theoremstyle{definition}
 
 \theoremstyle{remark}
 \newtheorem{rem}[lem]{Remark}
 \newtheorem{exa}[lem]{Example}
\begin{document}

\title{On the ramification of non-abelian
Galois coverings of degree $p^3$}

\author{Qizhi Zhang}
\maketitle

\begin{abstract}
The refined Swan conductor is defined by K.\ Kato \cite{KK2}, and
generalized by T.\ Saito \cite{wild}. In this part, we consider some
smooth $l$-adic \'{e}tale sheaves of rank $p$ such that we can be
define the $rsw$ following T.\ Saito, on some smooth dense open
subscheme $U$ of a smooth separated scheme X of finite type over a
perfect fields $\kappa$ of characteristic $p>0$. We give an explicit
expression of $rsw(\mathcal{F})$ in some situation. As a
consequence, we show that it is integral.
\end{abstract}

\section{Introduction}

%Where the Swan conductor $sw(\mathcal{F})$ measure the wild
%ramification of $\mathcal{F}$

The classical Swan conductor is defined in the case where the
extension of residue field is separable.

Kato in \cite{KK2}, gives a natural definition of the Swan conductor
for a character of degree one, which works without the separability
assumption of the extension of the residue field. In the same paper,
 he gives the definition of the refined Swan conductor $rsw$.

Kato in \cite{cft}, redefines the refined Swan conductor as an
element in the stalk of the sheaf $\omega _X ^1 (sw(\chi))$, where
$\omega _X ^1$ is the logarithmic cotangent sheaf, and extend it to
a global section. We call the property of rsw that can extend to a
global section \lq\lq \, integral " . He defines a $0$-cycle
$c_\chi$ to generalize the Swan conductor, by using the refined Swan
conductor.

 Along with the
generalization of the Swan conductor, Grothendieck-Ogg-Shafarevich
formula is generated also.

Kato and Saito in \cite{ramification}, define the Swan class as a
0-cycle class supported on the ramification locus as a refinement of
Swan conductor, and generalize the G-O-S formula to arbitrary rank.
In the situation that the scheme have dimension $\leq 2$ and the
sheave have rank 1, the Swan class is consistent with the $0$-cycle
class $c_\chi$ defined in \cite{cft}.

 Abbes and Saito in \cite{characteristic}, define a characteristic class of an $l$-adic sheaf as
 a refinement of the Euler-Poincar\'{e} characteristic, and get a
 refinement of the G-O-S formula for $l$-adic sheaf of rank $1$, in term characteristic class and
 the 0-cycle class defined for rank 1 sheaves in \cite{cft}, by explain the refine Swan conductor by a linear form on
some vector bundle on $D$.

Saito in \cite{wild}, by redefine rsw, defines the characteristic
cycle of an $l$-adic sheaf, as a cycle on the logarithmic cotangent
bundle. Under some condition (\cite{wild}, p53, the condition (R)
and (C)), he proves that the intersection of characteristic cycle
with the $0$-section computes the characteristic class defined in
\cite{characteristic} for $l$-adic sheaf of arbitrary rank, and
hence the Euler number. The property \lq\lq \, integral " of $rsw$
 is a part of the condition.

In this paper, we consider some smooth $l$-adic \'{e}tale sheaves of
rank $p$, that can be defined the $rsw$ follow \cite{wild}, on some
smooth dense open subscheme $U$ of a smooth separated scheme X of
finite type over a perfect fields $\kappa$ of characteristic $p>0$.
We give a explicit expression of $rsw(\mathcal{F})$ in some
situation. In fact, we would see, that it is integral, namely
\begin{displaymath}
rsw(\mathcal{F}) \in \Gamma (D, i_D ^\star \Omega ^1 _{\Delta X}
(log D)(kD))
\end{displaymath}
where $i_D : D \longrightarrow X$ is the closed immersion.

In \S 2, we study the ramification of the Artin-Schreier covering.
In \S 3, we study the ramification of some Galois coverings of
degree $p^3$. In \S 4, we study the ramification of some $l$-adic
sheaves.

\section{Artin-Schreier covering}

   Let $X$ be a separated smooth scheme of finite type over a perfect field
$\kappa$ of characteristic $p>0$, $D=X \setminus U$ be a smooth
divisor, $U=X \setminus D$.

Let $(X \times X)^\sim$ and $(X \times X)^{(kD)}$ be the log blow-up
and the log diagram blow-up of $X \times X$ defined in \cite{wild}
respectively, $E^0:=(X \times X)^\sim \setminus U \times U$, $E:=((X
\times X)^{(kD)} \setminus U \times U)_{red}$. We denote the ideal
sheaves of $\Delta X$ in $(X \times X)^\sim$ and $(X \times
X)^{(kD)}$ by $\mathcal{I}_{\Delta X}$ and $\mathcal{J}_{\Delta X}$
respectively. Let $\xi$, $\eta$ be the generic point of $D$ and $E$
respectively, then there exists only valuation $v_E$ of
$\mathcal{O}_{(X \times X)^{(kD)}, \eta}$, such that $v_E(f\otimes
1)=v_D(f)$ for any $f \in \mathcal{O}_{X, \xi}$.

\begin{lem}
\label{valuation}
% For $f \in \mathcal{O}_{U,\xi}$, if $v_D (f) =-n$, then
%as the element of $\mathcal{K}((X \times X)^{(kD)})$, $v_E (1
%\otimes f -f \otimes 1) \geq k-n$.
If $f \in \Gamma(X, \mathcal{O}(nD))$, then $1\otimes f -f \otimes 1
\in \Gamma((X \times X)^{(kD)}, \mathcal{J}_{\Delta X}((n-k)E))$.
\end{lem}
\noindent $\mathbf{Proof.}$ Let us consider the commutative diagram
\begin{displaymath}
\xymatrix@!=1.5pc{
                  E \ar[rr] \ar[d] && (X \times X)^{(kD)} \ar[d]^{\phi ^k} &&
                  \Delta X \cup kE \ar[d] \ar[ll] \\
                  E^0 \ar[rr] \ar[d] && (X \times X)^\sim \ar[d] &&
                  \Delta X \ar[ll] \\
                  D \ar[rr] && X \\
}
\end{displaymath}
We have
\begin{displaymath}
\mathcal{I}_{\Delta X} \longrightarrow \phi ^k _{\star}
\mathcal{J}_{\Delta X} (-kE), \quad \mathcal{O}_{(X \times X)^\sim}
(nE^0) \longrightarrow \phi ^k _{\star} \mathcal{O}_{(X \times
X)^{(kD)}}(nE)
\end{displaymath}
hence
\begin{displaymath}
\mathcal{I}_{\Delta X}(nE^0) \longrightarrow \phi ^k _{\star}
\mathcal{J}_{\Delta X} ((n-k)E)
\end{displaymath}

Therefore
\begin{displaymath}
\begin{array}{l}
1 \otimes f-f \otimes 1 = f \otimes 1 (f^{-1}\otimes f -1) \in \Gamma((X \times X)^\sim, \mathcal{I}_{\Delta X}(nE^0))
\\
\longrightarrow \Gamma((X \times X)^\sim, \phi ^k _\star
\mathcal{J}_{\Delta X}(n-k)E) \longrightarrow \Gamma((X \times
X)^{(kD)}, \mathcal{J}_{\Delta X}(n-k)E)
\end{array}
\end{displaymath}
\qed
\[\]
We have an exact sequence
\begin{displaymath}
H^0(U,\mathcal{O}_U) \longrightarrow H^0(U,\mathcal{O}_U)
\longrightarrow H^1(U, \mathbb{Z}/p\mathbb{Z}) \longrightarrow 0
\end{displaymath}
Therefore any Galois covering of $U$ with Galois group
$\mathbb{Z}/p\mathbb{Z}$ can be defined by equation $T^p-T-f$, where
$f \in H^0(U, \mathcal{O}_U)$. $T^p-T-f$ and $T^p-T-f_0$ define same
covering if and only if $f=f_0+f_1-f_1^p$ for some $f_1 \in H^0(U,
\mathcal{O}_U)$. It easy to see, for any Galois covering $V$ of $U$
with Galois group $\mathbb{Z}/p\mathbb{Z}$, that wild ramified on
$D$, there exists a $n \in \mathbb{N}$ such, that V can be defined
by equation $T^p-T-f$ for some $f \in \Gamma(X, \mathcal{O}_X(nD))$
and $df \in \Omega ^1 _{\Delta X} (log D)(nD) _\xi \setminus \Omega
^1 _{\Delta X} (log D)(nD^-)_\xi$.

\begin{prop}
\label{valuation}
 If $f \in \Gamma(X, \mathcal{O}(nD))$ and $df \in \Omega ^1 _{\Delta X} (log D)(nD) _\xi
\setminus \Omega ^1 _{\Delta X} (log D)(nD^-)_\xi$, then as the
element of $\mathcal{K}((X \times X)^{(R)})$, $v_E (1 \otimes f -f
\otimes 1) = k-n$.
\end{prop}

\noindent $\mathbf{Proof.}$ Consider the commutative diagram
\begin{displaymath} \xymatrix@!=1.5pc{
                  (X \times X)^{(kD)} && \Delta X \ar[ll]_(.4)\delta  \\
                  E \ar[u]^i && D \ar[ll]^\theta \ar[u]_{i_D} \\
                   && \xi \ar[u]
}
\end{displaymath}
Then we have a commutative diagram
\begin{displaymath}
\xymatrix{
                  \Gamma ((X \times X)^{(kD)}, \mathcal{J}_{\Delta X}((n-k)E)) \ar[rr] \ar[d] && \Gamma (\Delta X, \delta ^\star \mathcal{J}_{\Delta X}((n-k)E) ) \ar[d]
                  \\
                  \Gamma (E, i^\star \mathcal{J}_{\Delta X}((n-k)E) ) \ar[rr] && \Gamma (D, i_D ^\star \delta ^\star \mathcal{J}_{\Delta X}((n-k)E)
                  )
}
\end{displaymath}
but
\begin{displaymath}
\begin{array}{l}
\delta ^\star \mathcal{J}_{\Delta X}((n-k)E)=\delta ^\star
\mathcal{J}_{\Delta X} \cdot \delta ^\star \mathcal{O}_{(X \times
X)^{(kD)}}((n-k)E)=\delta ^\star \mathcal{J}_{\Delta X} \cdot
\mathcal{O}_{\Delta X}((n-k)D) \\
\longrightarrow \Omega _{\Delta X} ^1 (\log D)(kD) \cdot
\mathcal{O}_{\Delta X}((n-k)D)=\Omega _{\Delta X} ^1 (\log D)(nD)
\end{array}
\end{displaymath}
and
\begin{displaymath}
i^\star \mathcal{J}_{\Delta X}((n-k)E)=\mathcal{I}_D((n-k)E)
\end{displaymath}
Therefore we have
\begin{displaymath}
\xymatrix{
                  \Gamma ((X \times X)^{(kD)}, \mathcal{J}_{\Delta X}((n-k)E)) \ar[rr] \ar[d] && \Gamma (\Delta X, \delta ^\star \mathcal{J}_{\Delta X}((n-k)E) )
                  \ar[d] \ar[rr] && \Gamma(\Delta X, \Omega _{\Delta X} ^1 (\log
                  D)(nD)) \ar[d]
                  \\
                  \Gamma (E,  \mathcal{I}_{D}((n-k)E) ) \ar[rr] && \Gamma (D, i_D ^\star \delta ^\star \mathcal{J}_{\Delta
                  X}((n-k)E) \ar[rr] && \Gamma(D, i_D ^\star \Omega _{\Delta X} ^1 (\log D)(nD))
                  )
}
\end{displaymath}
The image of $1\otimes f-f\otimes 1$ in  $\Gamma(D, i_D ^\star
\Omega _{\Delta X} ^1 (\log D)(nD))$ is $df$.

However
\begin{displaymath}
\begin{array}{rl}
& i_D ^\star \Omega ^1 _{\Delta X} (log D)(kD) \\
= & i_D ^{-1} \Omega ^1 _{\Delta X} (log D)(kD) \otimes _{i_D ^{-1}
\mathcal{O} _{\Delta X}} \mathcal{O}_D \\
= & i_D ^{-1} \Omega ^1 _{\Delta X} (log D)(kD) \otimes _{i_D ^{-1}
\mathcal{O} _{\Delta X}} i^{-1} _D \frac{\mathcal{O}_{\Delta
X}}{\mathcal{I}_D} \\
= & \frac{i_D ^{-1} \Omega ^1 _{\Delta X} (log D)(kD)}{ i_D
^{-1}(\mathcal{I}_D \Omega ^1 _{\Delta X} (log D)(kD))} \\
= & i_D ^{-1} \frac{\Omega ^1 _{\Delta X} (log D)(kD)}{\Omega ^1
_{\Delta X} (log D)(kD^-)}
\end{array}
\end{displaymath}
hence
\begin{displaymath}
\begin{array}{rl}
&(i_D ^\star \Omega ^1 _{\Delta X} (log D)(kD))_\xi \\
=& (i_D ^{-1} \frac{\Omega ^1 _{\Delta X} (log D)(kD)}{\Omega ^1
_{\Delta X} (log D)(kD^-)})_\xi \\
= &(\frac{\Omega ^1 _{\Delta X} (log D)(kD)}{\Omega ^1 _{\Delta X}
(log D)(kD^-)})_\xi
\end{array}
\end{displaymath}
$df$ is not $0$ in $(i_D ^\star \Omega ^1 _{\Delta X} (log
D)(kD))_\xi$, so is not $0$ in $\Gamma(D, i_D ^\star \Omega _{\Delta
X} ^1 (\log D)(nD))$. Hence $1\otimes f-f \otimes 1$ is not $0$ in
$\Gamma (E, \mathcal{I}_{D}((n-k)E) )$. Therefore $v_E(1\otimes f-f
\otimes 1)=k-n$.

%Take a prime element $\pi$ of $\mathcal{O}_{X,\xi}$, then
%\begin{displaymath}
%\pi ^{n-k} df \in \Omega ^1 _{\Delta X} (log D)(kD) _\xi \setminus
%\Omega ^1 _{\Delta X} (log D)(kD^-)_\xi
%\end{displaymath}
%hence $\pi ^{k-n}df$ is not $0$ in $(\frac{\Omega ^1 _{\Delta X}
%(log D)(kD)}{\Omega ^1 _{\Delta X} (log D)(kD^-)})_\xi$. But this is
%the image of $\pi ^{n-k} \otimes 1 (f\otimes 1 -1 \otimes f)$ in
%$(\frac{\Omega ^1 _{\Delta X} (log D)(kD)}{\Omega ^1 _{\Delta X}
%(log D)(kD^-)})_\xi$. Hence $\pi ^{n-k} \otimes 1 (1\otimes f -f
%\otimes 1)$ is not $0$ in $\Gamma (E,\mathcal{I}_D)$ by lemma
%\eqref{diagrameomega}. So is not in $\Gamma (E,\mathcal{O}_E)$ also.
%Hence it is invertible  in $\mathcal{O}_{(X \times X)^{(kD)},\eta}$.
%Therefore
%\begin{displaymath}
%v_E (\pi ^{n-k} \otimes 1 (1\otimes f -f \otimes 1)) =0
%\end{displaymath}
%hence
%\begin{displaymath}
%v_E (1\otimes f -f \otimes 1) =k-n
%\end{displaymath}

%\begin{cor} \label{valuation} If $T^p-T-f_0=0$ is a minimal
%equation of the Galois covering $V=Spec
%\frac{\mathcal{O}_U[T]}{(T^p-T-f_0)}$, and $v_D(f_0)=-n$, then as
%the element of $\mathcal{K}((X \times X)^{(R)})$, $v_E (1 \otimes f
%-f \otimes 1) = k-n$. $\end{cor}

%\noindent $\mathbf{Proof.}$ Because lemma \ref{MinToOmega} and
%proposition \ref{diagrameomega}.

\section{Galois coverings of degree $p^3$}

Now, let us consider the Galois covering $V$ defined by equations
$T^p-T-f, S^p-S-g, U^p-U-fS-h$ on $U$, where $f, g, h \in
 \Gamma(U, \mathcal{O}_U)$. The Galois group $G=\left\{
                       \left(
                              \begin{array}{ccc}
                              1 & a & c\\
                              0 & 1 & b\\
                              0 & 0 & 1
                              \end{array}
                       \right)|a, b, c \in \mathbb{Z}/p\mathbb{Z}
               \right\}$.
$V$ has two quotient coverings $V_T$ and $V_S$, which are defined on
$U$ by equation $T^p-T-f$ and $S^p-S-g$ respectively. In the
following part of this paper, we always suppose
\begin{displaymath}
df \in \Omega ^1 _{X} (\log D) (nD) _\xi \setminus \Omega ^1 _{X}
(\log D) (nD^-) _\xi
\end{displaymath}
\begin{displaymath}
 dg \in \Omega ^1 _{X} (\log D) (mD) _\xi
\setminus \Omega ^1 _{X} (\log D) (mD^-) _\xi
\end{displaymath}
\begin{displaymath}
 dh \in \Omega ^1 _{X} (\log D) (rD) _\xi
\setminus \Omega ^1 _{X} (\log D) (rD^-) _\xi
\end{displaymath}
where $n=-v_D (f), m=-v_D(g), r=-v_D(h)$.

We have a filter of subgroups of $G \times G$ as follows :
\begin{displaymath}
G \times G \supset N \supset \Delta G \supset {1}
\end{displaymath}
where $N=\left\{
               \left(
                     \left(
                           \begin{array}{ccc}
                            1 & a & c \\
                             0 & 1 & b \\
                            0 & 0 & 1
                            \end{array}
                     \right),
                     \left(
                     \begin{array}{ccc}
                      1 & a & c' \\
                       0 & 1 & b \\
                       0 & 0 & 1
                     \end{array}
                     \right)
               \right)
               | a, b, c, c'\in\mathbb{F}_p
         \right\}
$, $\Delta G= \{ (g,g)| g \in G  \}$.

Let $Z_1:=(V \times V) _{\Delta G}, W_1:= (V \times V) _N$ be the
quotient of $V \times V$ under the action of $\Delta G$ and $N$
respectively; $Z_0, W_0$ be the normalization of $Z_1, W_1$ over $(X
\times X) ^{(kD)}$ respectively.

\begin{lem}

(1). If $k \geq \max \{m, n\}$, then $W_0 \longrightarrow (X \times
X)^{(kD)}$ is Galois with Galois group isomorphic to
$\mathbb{Z}/p\mathbb{Z} \times \mathbb{Z}/p\mathbb{Z} $

(2). If $k>\max \{m, n\}$, then $Z_0$ and $W_0$ is splitting to
$p^2$'s connected components over $E\cup \Delta X$.
\end{lem}

\noindent $\mathbf{Proof.}$

(1).
\begin{displaymath}
\begin{array}{rl}
V \times V=  & Spec \frac{\mathcal{O}_U \otimes \mathcal{O}_U
[T,S,U,T',S',U']}{(
\begin{array}{c}
T^p -T - f \otimes 1, S^p -S -g
\otimes 1, U^p -U-(f\otimes 1 S+h \otimes 1),\\
{T'}^p -T' -1\otimes f, {S'}^p -S' -1 \otimes g,{U'}^p -U' -(1
\otimes f S'+1 \otimes h)
\end{array}
)} \\
= & Spec \frac{\mathcal{O}_U \otimes \mathcal{O}_U
[T,S,U,\check{T},\check{S},\check{U}]}{(
\begin{array}{c}
T^p -T - f \otimes 1, S^p -S -g
\otimes 1, U^p -U-(f\otimes 1 S+h),\\
\check{T}^p -\check{T} -K_1, \check{S}^p -\check{S} -K_2,
\check{U}^p -\check{U} -(1 \otimes f S'-f \otimes 1 S +K_4)
\end{array}
)}
\end{array}
\end{displaymath}
Where $K_1 = 1 \otimes f -f\otimes 1, K_2= 1 \otimes g -g\otimes 1,
K_4= 1 \otimes h -h \otimes 1 \in \mathcal{O} _{U \times U}$

Therefore
\begin{displaymath} W_1= Spec \frac{\mathcal{O}_U \otimes
\mathcal{O}_U [\check{T},\check{S}]}{( \check{T}^p -\check{T} -K_1,
\check{S}^p -\check{S} -K_2 )}
\end{displaymath}

If $k \geq \max (m,n)$, then by lemma \eqref{valuation},  we know
$K_1, K_2 \in \Gamma ((X \times X)^{(kD)}, \mathcal{O}_{(X \times
X)^{(kD)}} )$. Then
\begin{displaymath}
W_0= Spec \frac{\mathcal{O}_{(X \times X)^{(R)}}
[\check{T},\check{S}]}{( \check{T}^p -\check{T} -K_1, \check{S}^p
-\check{S} -K_2 )}
\end{displaymath}

(2). If $k > \max (m,n)$, then by lemma \eqref{valuation}, we know
$v_E(K_1)>0, v_E(K_2)>0$. On the other hand, we know $K_{1|\Delta
X}=0$, $K_{2|\Delta X}=0$ because $K_{1|\Delta U}=0$, $K_{2|\Delta
U}=0$. Then
\begin{displaymath}
\begin{array}{rl}
W_0 \times _{(X \times X)^{(kD)}} (E\cup \Delta X )= & Spec
\frac{\mathcal{O}_{(X \times X)^{(R)}} [\check{T},\check{S}]}{(
\check{T}^p -\check{T} -K_1, \check{S}^p -\check{S} -K_2 )} \otimes
_{\mathcal{O}_{(X\times
X)^{(R)}}} \mathcal{O}_{(E\cup \Delta X)} \\
= & Spec \frac{\mathcal{O}_{(E\cup \Delta X)}
[\check{T},\check{S}]}{( \check{T}^p
-\check{T}, \check{S}^p -\check{S} )} \\
= & \coprod _{\check{t}, \check{s}} (E\cup \Delta
X)_{\check{t},\check{s}}
\end{array}
\end{displaymath}
Therefore
\begin{displaymath}
Z_0 \times _{(X \times X)^{(kD)}} (E\cup \Delta X )=Z_0 \times
_{W_0} \coprod _{\check{t}, \check{s}} (E\cup \Delta
X)_{\check{t},\check{s}} =\coprod _{\check{t}, \check{s}}
\tilde{F}_{\check{t},\check{s}}
\end{displaymath}

\qed

\begin{thm}
\label{L2d} Take a suitable $k$, then we have
%a commutative diagram
%\begin{align}
%\label{2d} \xymatrix@!=1.5pc{
%  Z_1 \ar[rr] \ar[dd]  && Z_0 \ar[dd] && \coprod _{\check{t},\check{s}} \tilde{F}_{\check{t},\check{s}} \ar[ll] \ar[dd] && \tilde{F}_{0,0} \ar[ll] \ar[dd]  && \coprod _r \Delta X_{0,0,r} \ar[ll] \ar[dd] && \Delta X_{0,0,0} \ar[ll] \ar[dd]\\
%  &&&&&&& F_{0,0} \ar[lu] \ar[dd] \\
%  W_1 \ar[rr] \ar[dd]_{\phi _1}  && W_0 \ar[dd] && \coprod _{\check{t},\check{s}} (E \cup \Delta X)_{\check{t},\check{s}} \ar[ll] \ar[dd] && (E \cup \Delta X)_{0,0} \ar[ll] \ar[dd] && \Delta X_{0,0} \ar'[l][ll] \ar[dd] && \Delta X_{0,0} \ar[ll] \\
% &&&&&&& E_{0,0} \ar[lu] \ar[dd]\\
% U \times U \ar[rr] ^j && (X \times X)^{(R)} && E\cup \Delta X \ar[ll] _{i^+} && (E\cup \Delta X )  \ar[ll] && \Delta X \ar'[l][ll] \\
% &&&&&&& E \ar[lu]
% }
%\end{align}
%satisfying the follow conditions

(1). Over $E_{0,0}$, $\tilde{F}$ is a Galois covering defined by an
Artin-Schreier equation, or splits to $p$'s connected components;
over $\Delta X_{0,0}$, it splits $p$'s connected components.

(2). In the commutative diagram
\begin{displaymath}
 \xymatrix@!=1.5pc{
  Z_1 \ar[rr] \ar[d]  && Z_0 \ar[d] && \coprod _{\check{t},\check{s}} \tilde{F}_{\check{t},\check{s}} \ar[ll] \ar[d] && \tilde{F}_{0,0} \ar[ll] \ar[d]  \\
  W_1 \ar[rr] \ar[d]_{\phi _1}  && W_0 \ar[d] && \coprod _{\check{t},\check{s}} (E \cup \Delta X)_{\check{t},\check{s}} \ar[ll] \ar[d] && (E \cup \Delta X)_{0,0} \ar[ll] \ar[d] \\
 U \times U \ar[rr] ^j && (X \times X)^{(R)} && E\cup \Delta X \ar[ll] _{i^+} && (E\cup \Delta X )  \ar[ll]
 }
\end{displaymath}
the base change map $\psi^{32\star} _{31} \psi ^{31\star} _{11} \psi
^{01} _{11\star} \longrightarrow  \psi ^{32\star} _{12} \psi ^{02}
_{12\star} \psi^{02\star} _{01} $ is a isomorphism for constructible
sheaf. Where we use $\psi ^{ij} _{kl}$ to denote the unique morphism
from the object at site $(i,j)$ to the object at site $(k,l)$ (if it
exists). For example $\psi ^{00} _{10} =j$, $\psi ^{20} _{10}=i^+$.
\end{thm}

\noindent $\mathbf{Proof.}$

(1). We can see
\begin{displaymath}
Z_1 = Spec \frac{\mathcal{O}_{W_1}
[\check{U}-S\check{T}]}{((\check{U}-S\check{T})^p
-(\check{U}-S\check{T})-K_3)}
\end{displaymath}
where $K_3=(K_1+f\otimes 1)\check{S}-g\otimes 1 (\check{T}+K_1) +
K_4$.

Let $\tilde{W}_0=W_0 \setminus \coprod _{(\check{t}, \check{s}) \neq
(0,0)} (E \cup \Delta X)_{\check{t},\check{s}}$,
$\tilde{Z}_0=\tilde{Z}_0 \times _{W_0} \tilde{W} _0 =Z_0 \setminus
\coprod _{(\check{t}, \check{s}) \neq (0,0)}
\tilde{F}_{\check{t},\check{s}}$. Then $\tilde{Z}_0$ is the
normalization of $Z_1$ over $\tilde{W}_0$. Let $\eta$ and $\eta
_{0,0}$ be the generic point of $E$ and $E_{0,0}$ respectively, then
there is a unique extension $v_{E_{0,0}}$ on the stalk
$\mathcal{O}_{\tilde{W}_0, \eta}$ of valuation $v_E$, such that the
restriction of $v_{E_{0,0}}$ on $\mathcal{O}_{(X \times X)^{(R)},
\eta}$ is just $v_E$.

We have
\begin{displaymath}
\begin{array}{rl}
K_3= & (K_1+f\otimes 1)(\check{S}^p-K_2)-g\otimes 1 (\check{T}+K_1)
+ K_4 \\
= & K_1 \check{S}^p + f \otimes 1 \check{S}^p -K_1K_2 -g\otimes 1
\check{T}^p +K_4 -f\otimes 1 K_2 \\
% = & K_3 ^+ + K_3 ^0
\end{array}
\end{displaymath}
where $v_{E_{0,0}}(f \otimes 1)=-n, v_{E_{0,0}}(g \otimes 1)=-m,
v_{E_{0,0}} (K_4)\geq k-r$ by lemma \ref{valuation}, and
$v_{E_{0,0}}(K_1)=k-n, v_{E_{0,0}}(K_2)=k-m$ by proposition
\ref{valuation}. Moreover, we have $v_{E_{0,0}}(\check{T})=k-n$, and
$v_{E_{0,0}}(\check{S})=k-m$, because $\check{T}^p -\check{T}
-K_1=0$ and $\check{S}^p -\check{S} -K_2=0$. Therefore we can take
enough big k such, that $v_{E_{0,0}}(K_3) \geqslant 0$, then $K_3
\in \Gamma (\tilde{W}_0,\mathcal{O}_{\tilde{W} _0})$. Therefore
\begin{displaymath}
\tilde{Z}_0 = Spec \frac{\mathcal{O}_{\tilde{W}_0}
[\check{U}-S\check{T}]}{((\check{U}-S\check{T})^p
-(\check{U}-S\check{T})-K_3)}
\end{displaymath}
\begin{displaymath}
\label{2d} \xymatrix@!=1.5pc{
  Z_1 \ar[rr] \ar[dd] \ar[rrrd]^{\psi ^l } && Z_0 \ar@{.>}[dd] && \coprod _{\check{t},\check{s}} \tilde{F}_{\check{t},\check{s}} \ar[ll] \ar@{.>}[dd] && \tilde{F}_{0,0} \ar[ll] \ar@{.>}[dd] \ar[llld]^{\psi ^r} && \coprod _{r} \Delta X_{0,0,r} \ar[ll] \ar[dd]  \\
 &&& \tilde{Z}_0 \ar@{.>}[lu]_{\psi ^c}  \ar[dd]^(.25)\psi &&&& F_{0,0} \ar[lu] \ar[dd] \\
  W_1 \ar@{.>}[rr] \ar[dd] \ar[rrrd] ^{\psi _l} && W_0 \ar@{.>}[dd] && \coprod _{\check{t},\check{s}} (E \cup \Delta X)_{\check{t},\check{s}} \ar@{.>}[ll] \ar@{.>}[dd] && (E\cup \Delta X)_{0,0} \ar@{.>}[ll] \ar@{.>}[dd] \ar[llld]^{\psi _r} && \Delta X_{0,0} \ar'[l][ll] \ar[dd] \\
 &&& \tilde{W}_0 \ar@{.>}[lu]_{\psi _c} &&&& E_{0,0} \ar[lu] \ar[dd] \\
 U \times U \ar[rr] ^j && (X \times X)^{(R)} && (E\cup \Delta X) \ar[ll] _{i^+} && (E\cup \Delta X)  \ar[ll] && \Delta X \ar'[l][ll] \\
 &&&&&&& E \ar[lu]
 }
\end{displaymath}

We can see $\tilde{Z}_0 $ is a Galois covering of $\tilde{W}_0$ and
\begin{displaymath}
F_{0,0}=\tilde{Z}_0 \times _{\tilde{W}_0} E_{0,0}=  Spec
\frac{\mathcal{O}_{E_{0,0}} [R]}{(R^p -R-\overline{K}_3)}
% = Spec \frac{\mathcal{O}_{E_{0,0}} [R]}{(R^p -R-\overline{K_3 ^0})}
\end{displaymath}
is a Galois covering of $E_{0,0}$.

%(4). We consider the commutative diagram
%\begin{displaymath}
%\xymatrix@!=1.5pc{
% W_1 \ar[dd] \ar[rr] && \tilde{W} _0 \ar[dd]\\
% & \coprod _p \Delta U \ar[lu] \ar[dd] \ar[rr] && \coprod _{p-1} \Delta U \sqcup \Delta X \ar[lu] \ar[dd] && \Delta X \ar[ll] & D \ar[l]   \\
% U \times U \ar[rr] && (X \times X)^{(R)} \\
% & \Delta U \ar[rr] \ar[lu] && \Delta X \ar[lu]\\
%}
%\end{displaymath}
%$K_3$ vanish on $\coprod _p \Delta U$, so is on $\coprod _{p-1}
%\Delta U \sqcup \Delta X$, so is on $\Delta X$, so is on D.
%Therefore
%\begin{displaymath}
%F_{0,0} \times _{E_{0,0}} D = Spec \frac{\mathcal{O}_D
%%[R]}{(R^p-R)}= \coprod _r D
%\end{displaymath}

On the other hands, It easy to see $K_3$ vanish on $\Delta U
_{0,0}$, then on $\Delta X _{0,0}$. Therefore
\begin{displaymath}
\tilde{Z}_0 \times _{\tilde{W}_0} \Delta X_{0,0} =  Spec
\frac{\mathcal{O}_{\Delta X_{0,0}} [R]}{(R^p -R)} = \coprod _{r}
\Delta X_{0,0,r}
\end{displaymath}

(2). We have
\begin{displaymath}
\begin{array}{rl}
\psi^{32\star} _{31}  \psi ^{31\star} _{11} \psi ^{01} _{11\star} =
& \psi_{31} ^{32\star} \psi_{r}^\star \psi_{c} ^\star \psi_{c
\star}\psi _{l \star} \\
= & \psi_{31} ^{32\star} \psi_{r}^\star \psi _{l \star}  \quad \quad
(\psi _c \mbox{ is a open immersion}) \\
= & \psi ^{r\star} \psi ^\star \psi _{l \star} \\
= & \psi ^{r\star} \psi ^l _{\star} \phi ^{02 \star} _{01} \quad \quad (\mbox{smooth base change})\\
= & \psi ^{32\star} _{12} \psi ^{02} _{12\star} \psi^{02\star} _{01}
 \quad (\psi ^c  \mbox{ is a open immersion})
\end{array}
\end{displaymath}
\qed
\begin{cor}
\label{basechange} For constructible sheaf, The base change map
$\psi^{52\star} _{51} \psi ^{51\star} _{11} \psi ^{01} _{11\star}
\longrightarrow \psi ^{52\star} _{12} \psi ^{02} _{12\star}
\psi^{02\star} _{01} $ is a isomorphism.
\end{cor}
%\begin{prop}
%If $r \neq m+n$,  $F_{0,0} \longrightarrow E_{0,0}$ is a Galois
%covering defined by a Artin-Schreier equation.
%\end{prop}
\begin{align}
\label{2d} \xymatrix@!=1.5pc{
  Z_1 \ar[rr] \ar[d]  && Z_0 \ar[d] && \coprod _{\check{t},\check{s}} \tilde{F}_{\check{t},\check{s}} \ar[ll] \ar[d] && \tilde{F}_{0,0} \ar[ll] \ar[d] && \coprod_{r} \Delta X _{0,0,r} \ar[ll] \ar[d] && \Delta X _{0,0} \ar[ll] \ar[d] \\
  W_1 \ar[rr] \ar[d]_{\phi _1}  && W_0 \ar[d] && \coprod _{\check{t},\check{s}} (E \cup \Delta X)_{\check{t},\check{s}} \ar[ll] \ar[d] && (E \cup \Delta X)_{0,0} \ar[ll] \ar[d] && \Delta X _{0,0} \ar[ll] \ar[d] && \Delta X _{0,0} \ar[ll] \\
 U \times U \ar[rr] ^j && (X \times X)^{(R)} && E\cup \Delta X \ar[ll] _{i^+} && (E\cup \Delta X )
 \ar[ll] && \Delta X \ar[ll]
 }
\end{align}
\qed

\[\]
\begin{cor}
\label{bigdiagram}

Take a suitable k, we have a commutative diagram
\begin{align}
\label{3d}
 \xymatrix@!=1.5pc{
 V \times V \ar[dd] \ar[rd] && V \times _{U} V \ar[ll] \ar@{.>}[dd] \ar[rd] && \coprod_{\check{t},\check{s},\check{u}} V_{\check{t},\check{s},\check{u}} \ar[ll] \ar@{.>}[dd] \ar[rd] && \coprod_{\check{t},\check{s},\check{u}} V_{\check{t},\check{s},\check{u}} \ar[ll] \ar@{.>}[dd] \ar[rd]  && \Delta V \ar[ll] \ar@{.>}[dd] \ar[rd] \\
 & Y_0 \ar[dd] && Y _{E \cup \Delta X} \ar[dddd] \ar[ll] && Y_3 \ar[dd] \ar[ll] && Y_4 \ar[dd] \ar[ll] && \Delta \overline{V} \ar[dd] \ar[ll] \\
 Z_1 \ar[dd] \ar[rd] &&  \coprod_{\check{t},\check{s},\check{u}} U_{\check{t},\check{s},\check{u}} \ar@{.>}[ll] \ar@{.>}[dd] && \coprod_{\check{u}} U_{0,0,\check{u}} \ar@{.>}[ll] \ar@{.>}[dd] \ar@{.>}[rd] && \coprod_{\check{u}} U_{0,0,\check{u}} \ar@{.>}[ll] \ar@{.>}[dd] \ar@{.>}[rd] && U_{0,0,0} \ar@{.>}[ll] \ar@{.>}[dd] \ar@{.>}[rd] \\
 & Z_0 \ar[dd] &&&& \tilde{F}_{0,0} \ar'[ll][llll] \ar[dd] && \coprod _r \Delta X \ar[dd] \ar[ll] && \Delta X \ar[dd] \ar[ll]   \\
   W_1 \ar[rd] \ar[dd]_{\phi _1}
      &  & \coprod_{\check{t},\check{s}} U_{\check{t},\check{s}} \ar@{.>}[dd] \ar@{.>}[ll] \ar@{.>}[rd]    && U_{0,0} \ar@{.>}[ll] \ar@{.>}[dd] \ar@{.>}[rd] && U_{0,0} \ar@{.>}[ll] \ar@{.>}[dd] \ar@{.>}[rd] && U_{0,0} \ar@{.>}[ll] \ar@{.>}[dd] \ar@{.>}[rd]  \\
 & W_0 \ar[dd]
      &  & \coprod_{\check{t},\check{s}} (E \cup \Delta X)_{\check{t},\check{s}} \ar[ll]\ar[dd] && (E \cup \Delta X)_{0.0} \ar[ll] \ar[dd] && \Delta X \ar[ll] \ar[dd] && \Delta X \ar[dd] \ar[ll] \\
   U\times U \ar[rd]_j
      &  & \Delta U \ar@{.>}[rd] \ar@{.>}[ll]          && \Delta U \ar@{.>}[ll] \ar@{.>}[rd]   && \Delta U \ar@{.>}[ll] \ar@{.>}[rd]  && \Delta U \ar@{.>}[ll] \ar@{.>}[rd]^{j_U}\\
&  (X \times X)^{(R)}  &  & E\cup \Delta X \ar[ll] _{i+}    && (E
\cup \Delta X) \ar[ll] && \Delta X \ar[ll] && \Delta X \ar[ll] }
\end{align}
\end{cor}

\noindent $\mathbf{Denote.}$ Let $\phi ^{ijk} _{pqr}$ denote the
unique map from the object at site $(i, j, k)$ to the object at site
$(p, q, r)$ (if it exists) of the diagram \eqref{3d}
 . For example $\phi ^{000}
 _{100}=i^+$, $\phi ^{010} _{000}=j$, $\phi ^{011} _{000} =\phi
 _1$, and $\phi ^{410} _{400}=j_U$.

\section{$l$-adic sheaves of rank $p$}
$G$ have a normal subgroup $H=\left\{
                       \left(
                              \begin{array}{ccc}
                              1 & a & c\\
                              0 & 1 & 0\\
                              0 & 0 & 1
                              \end{array}
                       \right)|a, b, c \in \mathbb{Z}/p\mathbb{Z}
               \right\}$. For any character $\chi : \mathbb{Z}/p\mathbb{Z} \longrightarrow \overline{\mathbb{Q}}_l
               ^\star$, we define a character
\begin{displaymath}
\begin{array}{cccc}
\tilde{\chi}: &  H & \longrightarrow & \overline{\mathbb{Q}}^\star
_l
\\
& \left(
      \begin{array}{ccc}
      1 & a & c \\
      0 & 1 & 0 \\
      0 & 0 & 1
      \end{array}
      \right)            & \mapsto &  \chi (c)
\end{array}
\end{displaymath}
We denote the induced representation $Ind _H ^G \tilde{\chi}$ by
$\rho$. In fact, we can see
\begin{displaymath}
\begin{array}{rccc}
\rho : &G & \longrightarrow & GL(\overline{\mathbb{Q}}_l ^p) \\
& \left(
   \begin{array}{ccc}
   1 & a & c \\
   0 & 1 & b \\
   0 & 0 & 1
   \end{array}
   \right)              & \mapsto & \left(
                                          \begin{array}{cccccc}
                                         &&& \chi ((p-b)a+c) && \\
                                            &&&& \ddots & \\
                                          &&&&& \chi ((p-1)a+c) \\
                                          \chi (c) &&&&& \\
                                            & \ddots &&&& \\
                                            && \chi ((p-b-1)a+c) &&&
                                            \end{array}
                                    \right)
\end{array}
\end{displaymath}

This defines a locally constant $\overline{\mathbb{Q}}_l-$ sheaf
$\mathcal{F}$ of rank $p$ on $U$.

\begin{prop}
Under the setting of the diagram
\begin{displaymath}
 \xymatrix@!=1.5pc{
&&& D \ar[ld] \ar[rd] && \Delta U \ar[ld] _{j_u} \\
 U\times U \ar[rd]^j && E \ar[rd] && \Delta X \ar[ld] \\
   & (X \times X)^{(R)} && E\cup \Delta X \ar[ll]_{i^+} \\
}
\end{displaymath}
There exists a sheaf $\mathcal{L}$ satisfying the following
conditions

(1). $\mathcal{L}$ is embedding in $i^{+ \star} j_\star
\mathcal{H}$.

(2). $\mathcal{L} _{|E}$ is a constant sheaf or a locally constant
sheaf defined by an Artin-Schreier equation which constant term a
linear form on $E$.

(3). $\mathcal{L} _{|\Delta X} = j_{U \star} \overline{\mathbb{Q}_l}
id $.
\end{prop}

\noindent $\mathbf{Proof.}$

(1). $\mathcal{H}$ is determined by the action of $\pi _1 (U \times
U)$ on $Hom(\overline{\mathbb{Q}}_l ^p, \overline{\mathbb{Q}}_l ^p)$
as follow:
\begin{displaymath}
\begin{array}{ccccc}
 \pi _1 (U \times U) & \longrightarrow & Gal(V \times V / U \times
 U) & \longrightarrow &
GL(Hom(\overline{\mathbb{Q}}_l ^p, \overline{\mathbb{Q}}_l ^p)) \\
&&
   \left(
   \left(
   \begin{array}{ccc}
   1 & a & c \\
   0 & 1 & b \\
   0 & 0 & 1
   \end{array}
   \right)
   ,
   \left(
   \begin{array}{ccc}
   1 & a' & c' \\
   0 & 1 & b' \\
   0 & 0 & 1
   \end{array}
   \right)
   \right)                      & \mapsto &  P' \times \quad
   \times P^{-1}
\end{array}
\end{displaymath}

where
\begin{displaymath}
P= \left(
          \begin{array}{cccccc}
          &&& \chi ((p-b)a+c) && \\
          &&&& \ddots & \\
          &&&&& \chi ((p-1)a+c) \\
          \chi (c) &&&&& \\
          & \ddots &&&& \\
          && \chi ((p-b-1)a+c) &&&
          \end{array}
  \right)
\end{displaymath}
\begin{displaymath}
P'= \left(
          \begin{array}{cccccc}
          &&& \chi ((p-b')a'+c') && \\
          &&&& \ddots & \\
          &&&&& \chi ((p-1)a'+c') \\
          \chi (c') &&&&& \\
          & \ddots &&&& \\
          && \chi ((p-b'-1)a'+c') &&&
          \end{array}
   \right)
\end{displaymath}.

For any $\sigma \in \pi _1 (W_1)$, the image of $\sigma$ in $Gal(V
\times V /U \times U)$ is in fact contained in $Gal(V \times V /
W_1)$. Therefore we can write the image by $   \left(
   \left(
   \begin{array}{ccc}
   1 & a & c \\
   0 & 1 & b \\
   0 & 0 & 1
   \end{array}
   \right)
   ,
   \left(
   \begin{array}{ccc}
   1 & a & c' \\
   0 & 1 & b \\
   0 & 0 & 1
   \end{array}
   \right)
   \right)$.

We can see the image of the identity element $I \in
Hom(\overline{\mathbb{Q}}_l ^p, \overline{\mathbb{Q}}_l ^p)$ under
the action of $\sigma$ is $\chi (c' - c)I$. It follow that
$\overline{\mathbb{Q}}_l I$ is a $\overline{\mathbb{Q}}_l [\pi _1 (W
_1)]$-submodule of $Hom(\overline{\mathbb{Q}}_l ^p,
\overline{\mathbb{Q}}_l ^p)$. Therefore we got a smooth subsheaf
$\mathcal{D}$ of $\phi _1 ^\star \mathcal{H}$ which is trivialized
by $Z_1$.

Let functor $\phi ^{201\star} _{001} \phi ^{011} _{001\star} $ act
on $0 \rightarrow \mathcal{D} \rightarrow \phi _1 ^\star
\mathcal{H}$, we get
\begin{displaymath}
0 \rightarrow \phi ^{201\star} _{001} \phi ^{011} _{001\star}
\mathcal{D} \rightarrow \phi ^{201\star} _{001} \phi ^{011}
_{001\star} \phi _1 ^\star \mathcal{H}
\end{displaymath}
But we know the base change map $\phi ^{001\star} _{000} j_\star
\mathcal{H} \rightarrow \phi ^{011} _{001\star} \phi _1 ^\star
\mathcal{H}$ is a isomorphism, because $\phi ^{001} _{000}$ is
Galois, so is smooth. Therefore we have
\begin{displaymath}
\phi ^{201\star} _{001} \phi ^{011} _{001\star} \phi _1 ^\star
\mathcal{H} = \phi ^{201\star} _{001} \phi ^{001\star} _{000}
j_\star \mathcal{H} = \phi ^{201\star} _{200} \phi ^{200\star}
_{000} j_\star \mathcal{H} = i^{+\star} j_\star \mathcal{H}
\end{displaymath}
Let $\mathcal{L}=\phi ^{201\star} _{001} \phi ^{011} _{001\star}
\mathcal{D}$, then $\mathcal{L}$ is embedded to $i^{+\star} j_\star
\mathcal{H}$.

(2). Consider the diagram \eqref{2d}, We have
\begin{displaymath}
  \mathcal{L}_{|{E_{0,0}} |F_{0,0}} =   (\psi ^{31\star} _{11}
\psi ^{01} _{11\star} \mathcal{D} )_{|{E_{0,0}} |F_{0,0}}= (\psi
^{32\star} _{12} \psi ^{31\star} _{11} \psi ^{01} _{11\star}
\mathcal{D} )_{|F_{0,0}} = (\psi ^{32\star} _{12} \psi ^{02}
_{12\star} \psi^{02\star} _{01} \mathcal{D})_{|F_{0,0}} \quad \mbox{
(by theorem } \ref{L2d},  (2) )
\end{displaymath}
But $\psi^{02\star} _{01} \mathcal{D}$ is a constant sheaf, then
$\mathcal{L}_{|{E_{0,0}} |F_{0,0}}$ is a constant sheaf also.

(3). Refer to the diagram \eqref{3d}, we have
\begin{displaymath}
\phi ^{402\star} _{401} \phi ^{401\star} _{001} \phi ^{011} _{001
\star} \mathcal{D} = \phi^{402 \star} _{002} \phi^{012} _{002\star}
\phi ^{012\star} _{011} \mathcal{D}
\end{displaymath}
by corollary  \ref{basechange}. On the other hand, it easy to see,
\begin{displaymath}
 \phi^{402\star}_{401} \phi^{411}_{401\star} \phi ^{411 \star}
_{011} \mathcal{D} =
 \phi^{412}_{402\star}
\phi^{412\star}_{012} \phi ^{012\star} _{011} \mathcal{D}
\end{displaymath}
But $\phi ^{012\star} _{011} \mathcal{D}$ is a constant sheaf, then
\begin{displaymath}
\phi^{402 \star} _{002} \phi^{012} _{002\star} \phi ^{012\star}
_{011} \mathcal{D} =
 \phi^{412}_{402\star}
\phi^{412\star}_{012} \phi ^{012\star} _{011} \mathcal{D}
\end{displaymath}
Therefore
\begin{displaymath}
\phi ^{402\star} _{401} \phi ^{401\star} _{001} \phi ^{011} _{001
\star} \mathcal{D} = \phi^{402\star}_{401} \phi^{411}_{401\star}
\phi ^{411 \star} _{011} \mathcal{D}
\end{displaymath}
i.e
\begin{displaymath}
\mathcal{L} _{|\Delta X} = j_{U \star} \overline{\mathbb{Q}_l} id
\end{displaymath}
\qed

\begin{rem}
In \cite{wild}, the rsw is defined to a injective (\cite{wild},
Corollary 1.3.4):
\begin{displaymath}
rsw: Hom (Gr_{log} ^k G_K, \mathbb{F}_p ) \longrightarrow \Omega _F
^1 (log) \otimes _F \mathbf{m}_{\overline{K}} ^{(-k)} /
\mathbf{m}_{\overline{K}} ^{(-k)+}
\end{displaymath}
where $K$ is the henselization of the stalk of $\mathcal{O}_X$ at
the generic point of $D$, F is its residue fields, $\overline{K}$ is
a separable closure of $K$. $\mathbf{m}_{\overline{K}} ^{(-k)}=\{ a
\in \overline{K} |v_K(a) \geq -k\}$ and $\mathbf{m}_{\overline{K}}
^{(-k)+}=\{ a \in \overline{K} |v_K(a) > -k\}$. There is a
filtration $(G_{K,log} ^a)_{a \in \mathbb{Q} _{\geq 0}}$ of
$G_K=\mbox{Gal}(\overline{K}/K)$, and $Gr_{log} ^a G_K$ is the
graded pieces $(G_{K,log} ^a/G_{K,log} ^{a+} )$.

In our case, we have a morphism of filter
\begin{displaymath}
\xymatrix{ G_{K, \log}^0 \ar[d] & \supset & G_{K, \log}^k \ar[d] &
\supset &
G_{K, \log}^{k+} \ar[d] \\
   G      & \supset &  G^k      & \supset & {1}
 }
\end{displaymath}
and $Gr_{log} ^k G_K=G^k$, where
\begin{displaymath}
G^k=  \left\{    \left(
      \begin{array}{ccc}
      1 & 0 & c \\
      0 & 1 & 0 \\
      0 & 0 & 1
      \end{array}
      \right)               | c\in \mathbb{Z}/p\mathbb{Z} \right\}
      \simeq \mathbb{Z}/p\mathbb{Z}
\end{displaymath}.

 $\rho$ is a Galois representation of
dimension $p$, whose unramified on U. Its restriction on $G_{K,\log}
^r$ factors by
\begin{displaymath}
\begin{array}{rccc}
 & Gr_{\log} ^k G_K & \longrightarrow & GL(\overline{\mathbb{Q}}_l ^p) \\
& \left(
   \begin{array}{ccc}
   1 & 0 & c \\
   0 & 1 & 0 \\
   0 & 0 & 1
   \end{array}
   \right)              & \mapsto & \chi (c) I_p
\end{array}
\end{displaymath}

This is a direct sum of $p$'s character $\chi : G_{K,\log} ^r
\longrightarrow \overline{\mathbb{Q}}_l ^\star$. $\chi$ can be
regarded as a element of $Hom (Gr_{log} ^r G_K, \mathbb{F}_p )$,
therefore we can define the refined Swan conductor of $\mathcal{F}$
by $rsw(\chi)$.

Now, let us consider the commutative diagram following:
\begin{displaymath}
\xymatrix@!=5pc{ \tilde{W}_0 \ar[rd]^{\tilde{\psi}} \\
                 & (X \times X)^{(kD)} & \Delta X \ar[l] ^\delta \ar[ull]_{\tilde{\delta}} \\
                & E \ar[u]^i \ar[uul]^{\tilde{i}}   & D
                \ar[l]_\theta \ar[u]_{i_D}
             %      E \ar[ru]^{\tilde{i}} \ar[r]^i & (X \times X)^{(kD)} &
              %     \Delta X \ar[l]_\delta \ar[lu]_{\tilde{\delta}}
}
\end{displaymath}

where $\tilde{W}_0=W_0 \setminus \coprod_{(\check{t},\check{s})\neq
(0,0)} (E \cup \Delta X)_{(\check{t},\check{s})} $, and $\tilde{i},
\tilde{\delta}$ are the only liftings of $i, \delta$ respectively.
Let the ideal sheaf of $\Delta X$ in $\tilde{W}_0$ is
$\tilde{\mathcal{J}}_{\Delta X}$, then we have a commutative diagram
\begin{displaymath}
\xymatrix@!=3.25pc{
                  \Gamma (\tilde{W}, \tilde{\mathcal{J}}_{\Delta X}) \ar[rr] \ar[d] && \Gamma (\Delta X, \tilde{\delta} ^\star \tilde{\mathcal{J}}_{\Delta X} ) \ar[d]
                  \\
                  \Gamma (E, \tilde{i}^\star \tilde{\mathcal{J}}_{\Delta X} ) \ar[rr] && \Gamma (D, i_D ^\star \tilde{\delta} ^\star \tilde{\mathcal{J}}_{\Delta X}
                  ) \\
 }
\end{displaymath}
However, $\tilde{\mathcal{J}}_{\Delta X}=\tilde{\psi}^\star
\mathcal{J}_{\Delta X}$. Hence
\begin{displaymath}
\tilde{i}^\star \tilde{J}_{\Delta X} = i ^\star J_{\Delta X}, \quad
\tilde{\delta}^\star \tilde{J}_{\Delta X} = \delta ^\star J_{\Delta
X}
\end{displaymath}

Therefore, we have a commutative diagram
\begin{displaymath}
\xymatrix@!=3.25pc{
                  \Gamma (\tilde{W}, \tilde{\mathcal{J}}_{\Delta X}) \ar[rr] \ar[d] && \Gamma (\Delta X, \delta ^\star \mathcal{J}_{\Delta X} ) \ar[d] \ar[rr] && \Gamma (\Delta X, \Omega ^1 _{\Delta X} (log D)(kD) ) \ar[d] \\
                  \Gamma (E, \mathcal{I}_{D} ) \ar[rr] && \Gamma (D, i_D ^\star \delta ^\star \mathcal{J}_{\Delta X}
                  ) \ar[rr] && \Gamma (D, i_D ^\star \Omega ^1 _{\Delta X} (log D)(kD))  \\
                   % &&&& (\frac{\Omega ^1 _{\Delta X} (log D)(kD)}{\Omega ^1 _{\Delta X}(log D)(kD^-)})_\xi
 }
\end{displaymath}
By Corollary 1.3.4 and the proof of Theorem 1.3.3 in \cite{wild},
when $\mathcal{L}_{|E}$ is a locally constant sheaf defined by a
Artin-Schreier equation $R^p -R-\overline{K}_3$, the image of
$\overline{K}_3$ in $\Gamma (D, i_D ^\star \Omega ^1 _{\Delta X}
(log D)(kD))$ is the $rsw(\mathcal{F})$
\end{rem}

\begin{thm}
We can compute $rsw(\mathcal{F})$ as follow:

(1). If $r>m+n$, then $rsw(\mathcal{F})=dh \in \Gamma (D, i_D ^\star
\Omega ^1 _{\Delta X} (log D)(kD))$, where $k=r$.

(2). If $r<m+n$, then $rsw(\mathcal{F})=-fdg \in \Gamma (D, i_D
^\star \Omega ^1 _{\Delta X} (log D)(kD))$, where $k=m+n$.

(3). If $r=m+n$, and $n \geq m$, and $dh-fdg \in \Omega _X ^1 (log
D) (kD) _\xi \setminus \Omega _X ^1 (log D) (kD^-) _\xi$ for some
$k>n+\frac{m}{p}$, then $rsw(\mathcal{F})= dh-fdg \in \Gamma (D, i_D
^\star \Omega ^1 _{\Delta X} (log D)(kD))$.

(4). If $r=m+n$, and $m \geq n$, and $dh-fdg \in \Omega _X ^1 (log
D) (kD) _\xi \setminus \Omega _X ^1 (log D) (kD^-) _\xi$ for some
$k>m+\frac{n}{p}$, then $rsw(\mathcal{F})= dh-fdg \in \Gamma (D, i_D
^\star \Omega ^1 _{\Delta X} (log D)(kD))$.

\end{thm}

\noindent $\mathbf{Proof.}$We have
\begin{displaymath}
\begin{array}{rl}
K_3= & (K_1+f\otimes 1)(\check{S}^p-K_2)-g\otimes 1 (\check{T}+K_1)
+ K_4 \\
= & K_1 \check{S}^p + f \otimes 1 \check{S}^p -K_1K_2 -g\otimes 1
\check{T}^p +K_4 -f\otimes 1 K_2 \\
\end{array}
\end{displaymath}
where $v_{E_{0,0}}(f \otimes 1)=-n, v_{E_{0,0}}(g \otimes 1)=-m$,
and $v_{E_{0,0}} (K_4)= k-r, v_{E_{0,0}}(K_1)=k-n,
v_{E_{0,0}}(K_2)=k-m$ by proposition \ref{valuation}. Moreover, we
have $v_{E_{0,0}}(\check{T})=k-n$, and $v_{E_{0,0}}(\check{S})=k-m$,
because $\check{T}^p -\check{T} -K_1=0$ and $\check{S}^p -\check{S}
-K_2=0$. We can abuse $v_{E_{0,0}}$ and $v_E$, and write the
valuation of terms of $K_3$ in following table

\begin{tabular}{|c|c|c|c|c|c|c|}
\hline  term & $K_1 \check{S}^p$ & $f \otimes 1 \check{S}^p$ &
$K_1K_2$ & $g\otimes 1 \check{T}^p$ &  $K_4 $ & $f\otimes 1 K_2$ \\
\hline $v_E$ & $k-n+p(k-m)$ & $p(k-m)-n$ & $2k-m-n$ & $p(k-n)-m$ &
$k-r$ &
$k-m-n$ \\
\hline
\end{tabular}

Therefore,

(1). In this situation, $v_E(K_4)=0$, and $v_{E}(K_1 \check{S}^p + f
\otimes 1 \check{S}^p -K_1K_2 -g\otimes 1 \check{T}^p -f\otimes 1
K_2)>0$. Hence $\overline{K_3}=\overline{K_4} \in \Gamma(E,
\mathcal{I}_D)$. The image of $K_4$ in $\Gamma (D, i_D ^\star \Omega
^1 _{\Delta X} (log D)(kD))$ is $dh$.

(2). In this situation, $v_E(-f\otimes 1 K_2)=0$, and $v_E(K_1
\check{S}^p + f \otimes 1 \check{S}^p -K_1K_2 -g\otimes 1
\check{T}^p + K_4 )>0$. Hence $\overline{K_3}=\overline{-f\otimes 1
K_2} \in \Gamma(E, \mathcal{I}_D)$. The image of $-f\otimes 1 K_2$
in $\Gamma (D, i_D ^\star \Omega ^1 _{\Delta X} (log D)(kD))$ is
$-fdg$.

(3) By the following lemma, we know $v_E(K_4-f\otimes 1 K_2) \geq
0$. In other words, we know also
\begin{displaymath}
v_{E}(K_1 \check{S}^p +f \otimes 1 \check{S}^p -K_1K_2 -g\otimes 1
\check{T}^p )>0
\end{displaymath}
Hence $\overline{K_3}=\overline{K_4-f \otimes 1 K_2} \in \Gamma(E,
\mathcal{I}_D)$. The image of $K_4-f \otimes 1 K_2$ in $\Gamma (D,
i_D ^\star \Omega ^1 _{\Delta X} (log D)(kD))$ is $dh-fdg$. $dh-fdg$
is not $0$ as a element in $\Gamma (D, i_D ^\star \Omega ^1 _{\Delta
X} (log D)(kD))$, because it is not in $\Omega _X ^1 (log D) (kD^-)
_\xi$. Therefore $\overline {K_4}$ is not  $0$ in $\Gamma(E,
\mathcal{I}_D)$.

 (4). It is similar to (3).

\begin{lem}
Under the condition in (3) of previous theorem, we have
$v_E(K_4-f\otimes 1 K_2) \geq 0$.
\end{lem}
\noindent $\mathbf{Proof.}$ Let us consider the diagram
\begin{displaymath}
\xymatrix@!=3.25pc{
                   (X \times X)^{(kD)} \ar[d]^{\phi ^k} && \Delta X \cup kE \ar[d] \ar[ll] && kE \ar[ll] \\
                     (X \times X)^\sim && \Delta X
                     \ar[ll] _{\delta ^0}
               }
\end{displaymath}

%%%%%%%%%%%%%%%%%%%%%%%%%%%%%%%%%%%%%%%%%%%%%%%%%%%%%%%%%%%%%%%%%%%%%%%%%%%%%%%%%%%%%%%%%%%%%%%%%%%%%%%%%%%%%%%%%%%%%%%%%%%%%%%%%%%%%%%%%%%%%%%%%%%%%%%%%%%%%%%%%%%%%%%%%%%%%%%%%%%%%%%%%%%%%%%%%%%%%%%%%%%%%%%%%%%%%%%%%%%%%%%%%%%%%%%%%%%%%%%%%%%%%%%%%%%%%%%%%%%%%%%%%%%%%%%%%%%%%%%%%%%%%%%%%%%%%%%%%%%%%%%%%%%%%%%%%%%%%%%%%%%%%%%%%%%%%%%%%%%%%%%%%%%%%%%%%%%%%%%%%%%%%%%%%%%%%%%%%%%%%%%%%%%%%%%%%%%%%%%%%%%%%%%%%%%%%%%%%%%%%%%%%%%%%%%%%%%%%%%%%%%%%%%%%%%%%%%%55
Let $i^k$ be the morphism from $kE$ to $(X \times X)^{(kD)}$, then
we have a commutative diagram
\begin{displaymath}
\xymatrix@!=3pc{
                  0 \ar[rr]  && \mathcal{O}_{(X \times X)^{(kD)}}(-kE) \ar[d] \ar[rr] && \mathcal{O}_{(X \times X)^{(kD)}} \ar[rr] \ar[d] && i^k _\star \mathcal{O}_{kE} \ar[rr] \ar[d] && 0 \\
                  0 \ar[rr]  && \mathcal{O}_{(X \times X)^{(kD)}}((r-2k)E) \ar[rr] && \mathcal{O}_{(X \times X)^{(kD)}}((r-k)E) \ar[rr] && i^k _\star \mathcal{O}_{kE}((r-k)E) \ar[rr] && 0 \\
               }
\end{displaymath}
of sheaves of $\mathcal{O}_{(X \times X)^{(kD)}}$ module, and a
commutative diagram
\begin{displaymath}
\xymatrix@!=3pc{
                  0 \ar[rr]  && \mathcal{I}_{\Delta X} ^2 (kE^0) \ar[d] \ar[rr] && \mathcal{I}_{\Delta X} (kE^0) \ar[rr] \ar[d] && \delta ^0 _\star (\Omega _{\Delta X} ^1 (\log D)(kD)) \ar[rr] \ar[d] && 0 \\
                  0 \ar[rr]  && \mathcal{I}_{\Delta X} ^2 (rE^0)  \ar[rr] && \mathcal{I}_{\Delta X} (rE^0) \ar[rr]  && \delta ^0 _\star (\Omega _{\Delta X} ^1 (\log D)(rD)) \ar[rr]  && 0 \\
               }
\end{displaymath}
of sheaves of $\mathcal{O}_{(X \times X)^\sim}$ module. Where  where
$E^0$ is the complement of $U \times U$ in $(X \times X)^\sim$ as a
reduced scheme, and $\mathcal{I}_{\Delta X}$ is the ideal sheaf of
$\Delta X$ in $(X \times X)^\sim$. Let $\mathcal{J}_{\Delta X}$ be
the ideal sheaf of $\Delta X$ in $(X \times X)^{(kD)}$, then we have
a morphism
\begin{displaymath}
\mathcal{I}_{\Delta X} \longrightarrow \phi ^k _{\star}
\mathcal{J}_{\Delta X} (-kE) \longrightarrow \phi ^k _{\star}
\mathcal{O}_{(X \times X)^{(kD)}} (-kE)
\end{displaymath}
of sheaves of $\mathcal{O}_{(X \times X)^\sim}$ module. Therefore we
have a commutative disgrame
\begin{displaymath}
\xymatrix@!=3.1pc{
                  0 \ar[rr]  && \mathcal{I}_{\Delta X} ^2 (kE^0) \ar'[d][dd] \ar[rr] \ar[rd] && \mathcal{I}_{\Delta X} (kE^0) \ar[rr] \ar'[d][dd] \ar[rd] && \delta ^0 _\star (\Omega _{\Delta X} ^1 (\log D)(kD)) \ar[rr] \ar'[d][dd] \ar[rd] && 0 \\
                  &0 \ar[rr] && \phi ^k _\star \mathcal{O}_{(X \times X)^{(kD)}}(-kE) \ar[dd] \ar[rr]   && \phi ^k _\star \mathcal{O}_{(X \times X)^{(kD)}} \ar[rr] \ar[dd]  && \phi ^k _{\star} i^k _\star \mathcal{O}_{kE}  \ar[dd]   \\
                  0 \ar[rr]  && \mathcal{I}_{\Delta X} ^2 (rE^0)  \ar'[r][rr] \ar[rd] && \mathcal{I}_{\Delta X} (rE^0) \ar'[r][rr] \ar[rd] && \delta ^0 _\star (\Omega _{\Delta X} ^1 (\log D)(rD)) \ar'[r][rr] \ar[rd] && 0 \\
                  &0 \ar[rr] && \phi ^k _{\star} \mathcal{O}_{(X \times X)^{(kD)}}((r-2k)E) \ar[rr]  && \phi ^k _\star \mathcal{O}_{(X \times X)^{(kD)}}((r-k)E) \ar[rr]  && \phi ^k _\star i^k _\star \mathcal{O}_{kE}((r-k)E)    \\
               }
\end{displaymath}
of sheaves of $\mathcal{O}_{(X \times X)^\sim}$ module, where the
morphisms
\begin{displaymath}
\delta ^0 _\star (\Omega _{\Delta X} ^1 (\log D)(kD))
\longrightarrow \phi ^k _\star i^k _\star \mathcal{O}_{kE}, \quad
\delta ^0 _\star (\Omega _{\Delta X} ^1 (\log D)(rD))
\longrightarrow \phi ^k _\star i^k _\star \mathcal{O}_{kE}((r-k)E)
\end{displaymath}
are induced by
\begin{displaymath}
\mathcal{I}_{\Delta X} (kE^0) \longrightarrow \phi ^k _\star
\mathcal{O}_{(X \times X)^{(R)}} , \quad \mathcal{I}_{\Delta X}
(rE^0) \longrightarrow \phi ^k _\star \mathcal{O}_{(X \times
X)^{(R)}}((r-k)E)
\end{displaymath}
respectively.

Take the global sections, we have
\begin{displaymath}
\xymatrix@!=5pc{
                 % 0 \ar[rr]  && \Gamma((X \times X)^\sim, \mathcal{I}_{\Delta X} ^2 (kE^0)) \ar'[d][dd] \ar[rr] \ar[rd] &&
                  \Gamma((X \times X)^\sim, \mathcal{I}_{\Delta X} (kE^0)) \ar[rr] \ar'[d][dd] \ar[rd] && \Gamma(\Delta X, \Omega _{\Delta X} ^1 (\log D)(kD)) \ar'[d][dd] \ar[rd]\\
                  %0 \ar[rr] && \Gamma((X \times X)^{(kD)}, \mathcal{O}_{(X \times X)^{(kD)}}(-kE)) \ar[dd] \ar[rr]   &&
                  & \Gamma((X \times X)^{(kD)}, \mathcal{O}_{(X \times X)^{(kD)}}) \ar[rr] \ar[dd]  && \Gamma(kE, \mathcal{O}_{kE})  \ar[dd]   \\
                  %0 \ar[rr]  && \mathcal{I}_{\Delta X} ^2 (rE^0)  \ar'[r][rr] \ar[rd] &&
                  \Gamma((X \times X)^\sim, \mathcal{I}_{\Delta X} (rE^0)) \ar'[r][rr] \ar[rd] && \Gamma(\Delta X, \Omega _{\Delta X} ^1 (\log D)(rD)) \ar[rd]  \\
                  %&0 \ar[rr] && \phi ^k _{\star} \mathcal{O}_{(X \times X)^{(kD)}}((r-2k)E) \ar[rr]  &&
                  &\Gamma((X \times X)^{(kD)}, \mathcal{O}_{(X \times X)^{(kD)}}((r-k)E)) \ar[rr]  && \Gamma(kE, \mathcal{O}_{kE}((r-k)E))    \\
               }
\end{displaymath}

Now, $K_4$ and $f\otimes 1 K_2$ is in $\Gamma((X \times X)^\sim,
\mathcal{I}_{\Delta X} (rE^0))$, the image $dh-fdg$ of $K_4-f
\otimes 1 K_2$ in $\Gamma(\Delta X, \Omega _{\Delta X} ^1 (\log
D)(rD))$ is in fact inside $\Gamma(\Delta X, \Omega _{\Delta X} ^1
(\log D)(kD))$. Hence its image $\overline{K}_4-\overline{f\otimes 1
K_2}$ in $\Gamma(kE, \mathcal{O}_{kE}((r-k)E))$ is in fact inside
$\Gamma(kE, \mathcal{O}_{kE})$. Therefore $v_E(K_4-f\otimes 1 K_2)
\geq 0$.

\noindent $\mathbf{Acknowledgement.}$ I would like to show my
highest respect and appreciation to my advisor professor Takeshi
Saito for introducing this problem and giving some ideas of the
proof to me. I have been learning so much from him and without his
help, I could never have finished my master thesis.

%We need only to show, $\overline{K_3} \neq 0$ as a element in
%$\Gamma (E, \mathcal{O}_{E})$.

%By lemma \eqref{diagrameomega}, we can see $\overline{K_3 ^0}$ and
%$dh-fdg$ are the images of $K_3 ^0$ in $\Gamma (E, \mathcal{I}_D)$
%and $(\frac{\Omega ^1 _{\Delta X} (log D)(kD)}{\Omega ^1 _{\Delta X}
%(log D)((k-1)D)})_\xi $ respectively. Therefore $\overline{K_3 ^0}$
%is not $0$ in $\Gamma (E, \mathcal{I}_D)$, so is not also in $\Gamma
%(E, \mathcal{O}_E)$.

\end{document}